\documentclass{elsart}
\usepackage{amssymb}
\usepackage{graphicx}

\begin{document}
\begin{frontmatter}

\title{Parameter Collapse due to the Zeros in the Inverse Condition}
\author{R.~Spjut\thanksref{label2}}
\ead{spjut@math.ucsb.edu}
\ead[url]{http://math.ucsb.edu/$\sim$spjut/}
\thanks[label2]{Work supported by the UCSB and the NSF.  Special thanks to Mihai Putinar.}
\date{7 Feb 2008}
\address{Department of Mathematics, University of California, Santa Barbara, CA 93106}

\begin{abstract}
Helton, Lasserre, and Putinar (2008, Ann. Probability; arXiv:math/0702314) expose the relationship between three properties of a measure:
the \textit{conditional triangularity property} of the associated orthogonal polynomials,
the \textit{zeros in the inverse condition} of the truncated moment matrix,
and \textit{conditional independence}.
The purpose of this article is to provide examples of
\textit{parameter collapse to product structure} given that the zeros in the inverse condition holds up to some degree $d$.
  Specifically,
start with a parameterized family of probability density functions;
require that the \textit{zeros in the inverse condition} up to degree $d$ holds;
and validate that imposing this restriction on the parameterized family results in a measure with product structure,
or at least that conditional independence holds.
Algorithms related to parameter collapse are supplied, including the computation of the \textit{zeros in the inverse condition} up to degree $d$.
\end{abstract}

\end{frontmatter}

\section{Introduction}
Generalization of independence certificates for multivariate probability distributions is exciting.
For instance, techniques such as those found in \cite{2} could be generalized to non-normal variables.
As a specific example, near the variable bounds of instrument operation error distributions are known to be non-Gaussian.
In short, if a data set is fitted by a particular family of multivariate distributions,
the algorithms supplied in this article provide a statistician with the knowledge of precisely how many moments must be calculated
to ensure variables are independent, or at least conditionally independent.

This is perhaps the most interesting aspect,
that parameter collapse to product structure can require computation of moment matrices of order greater than 2.
To exhibit independence of variables for multivariate Gaussian distributions one much check the truncated moments only up to order 2,
but proving independence of variables for some distributions requires higher moments.
When a parameterized family of distributions exhibits parameter collapse,
the minimum order of the zeros in the inverse condition required for the measure to exhibit product structure is called
the \textit{parameter collapse order} of the parameterized measure \cite[\S6 Definition 3]{1}.

Finally, the computations of moments of polynomials by hand is tedious, at best.
Maple routines are provided in the appendix to automate computations.
Furthermore, a program is available upon request that computes the \textit{parameter collapse order} of a particular family of multivariate measures.
Assistance for translation to other computer algebra systems or languages is offered and code is immediately available in Maple, Mathematica, Matlab, C, and Java.
In particular, it may prove beneficial to integrate these routines with existing, high-powered software \cite{3}.

The following examples are meant to provide elementary exposition of the analysis discussed above.

\section{Examples}
\subsection{Preliminary example}

Consider the following parameterized multivariate distribution with product structure.  Let $X_{1},X_{2},X_{3}$ be independent Poisson processes on $\mathbb{R}_{+}=(0,\infty)$ with probability density
\[p_{\lambda,k}(x)=\left\{\begin{array}{cc}
                  \frac{\lambda^{k}x^{k-1}e^{-\lambda x}}{\Gamma(k)} & x\geq 0 \\
                  0 & x<0
                \end{array}\right.; \Gamma(k)=\int_{0}^{\infty} t^{k-1}e^{-t}dt.\]
Let $\lambda\tau=x$ and for simplicity, $\lambda =1$, so that the joint distribution is
\[p_{k_{1},k_{2},k_{3}}(x_{1},x_{2},x_{3})=\frac{e^{-x_{1}}x_{1}^{k_{1}-1}}{\Gamma(k_{1})}\frac{e^{-x_{2}}x_{2}^{k_{2}-1}}{\Gamma(k_{2})}\frac{e^{-x_{3}}x_{3}^{k_{3}-1}}{\Gamma(k_{3})},\]
a measure with product structure.  To compute moments, we are led to the gamma function
\begin{eqnarray*}
  E(x_{1}^{p}x_{2}^{q}) &=& \int_{0}^{\infty}\int_{0}^{\infty} x_{1}^{p}x_{2}^{q}\frac{e^{-x_{1}}x_{1}^{k_{1}-1}}{\Gamma(k_{1})}\frac{e^{-x_{2}}x_{2}^{k_{2}-1}}{\Gamma(k_{2})}dx_{1}dx_{2}\\
   &=& \int_{0}^{\infty}\frac{e^{-x_{1}}x_{1}^{p+k_{1}-1}}{\Gamma(k_{1})}dx_{1}\int_{0}^{\infty}\frac{e^{-x_{2}}x_{2}^{q+k_{2}-1}}{\Gamma(k_{2})}dx_{2}\\
   &=&\frac{\Gamma(p+k_{1})\Gamma(q+k_{2})}{\Gamma(k_{1})\Gamma(k_{2})}=p!q!.\\
\end{eqnarray*}
For $X_{1}$ and $X_{2}$, $M_{1}$ and $M_{1}^{-1}$ are:
\[  \begin{array}{ccc}\begin{array}{c|ccc}
            & 1 & x_{1} & x_{2} \\
      \hline
      1     & 1 & 1 & 1 \\
      x_{1} & 1 & 2 & 1 \\
      x_{2} & 1 & 1 & 2
    \end{array} &\hspace{20pt} &M_{1}^{-1} =
\left[ \begin {array}{ccc} 3&-1&-1\\\noalign{\medskip}-1&1&0\\\noalign{\medskip}-1&0&1\end {array} \right]\end{array}.\]
From \cite[Theorem 5]{1}, the zeros in the inverse property is equivalent to the orthogonal polynomials satisfying the conditional triangularity.
For $X_{1}$ and $X_{2}$, $M_{2}$  and $M_{2}^{-1}$ are:
\[  \begin{array}{ccc}\begin{array}{c|cccccc}
                 & 1 & x_{1} & x_{2} & x_{1}^{2} & x_{1}x_{2} & x_{2}^{2} \\
      \hline
      1          & 1 & 1     & 1     & 2         & 1          & 2 \\
      x_{1}      & 1 & 2     & 1     & 6         & 2          & 2 \\
      x_{2}      & 1 & 1     & 2     & 2         & 2          & 6 \\
      x_{1}^{2}  & 2 & 6     & 2     & 24        & 6          & 4 \\
      x_{1}x_{2} & 1 & 2     & 2     & 6         & 4          & 6 \\
      x_{2}^{2}  & 2 & 2     & 6     & 4         & 6          & 24
    \end{array}  &\hspace{20pt} &M_{2}^{-1} =
\left[ \begin {array}{cccccc} 6&-4&-4&1/2&1&1/2\\\noalign{\medskip}-4&6&1&-1&-1&0\\\noalign{\medskip}-4&1&6&0&-1&-1\\\noalign{\medskip}1/2&
-1&0&1/4&0&0\\\noalign{\medskip}1&-1&-1&0&1&0\\\noalign{\medskip}1/2&0
&-1&0&0&1/4\end {array} \right].\end{array}\]

\begin{figure}[htp]
\label{ZIIMM}
\centering
\includegraphics[width=142pt,height=116pt]{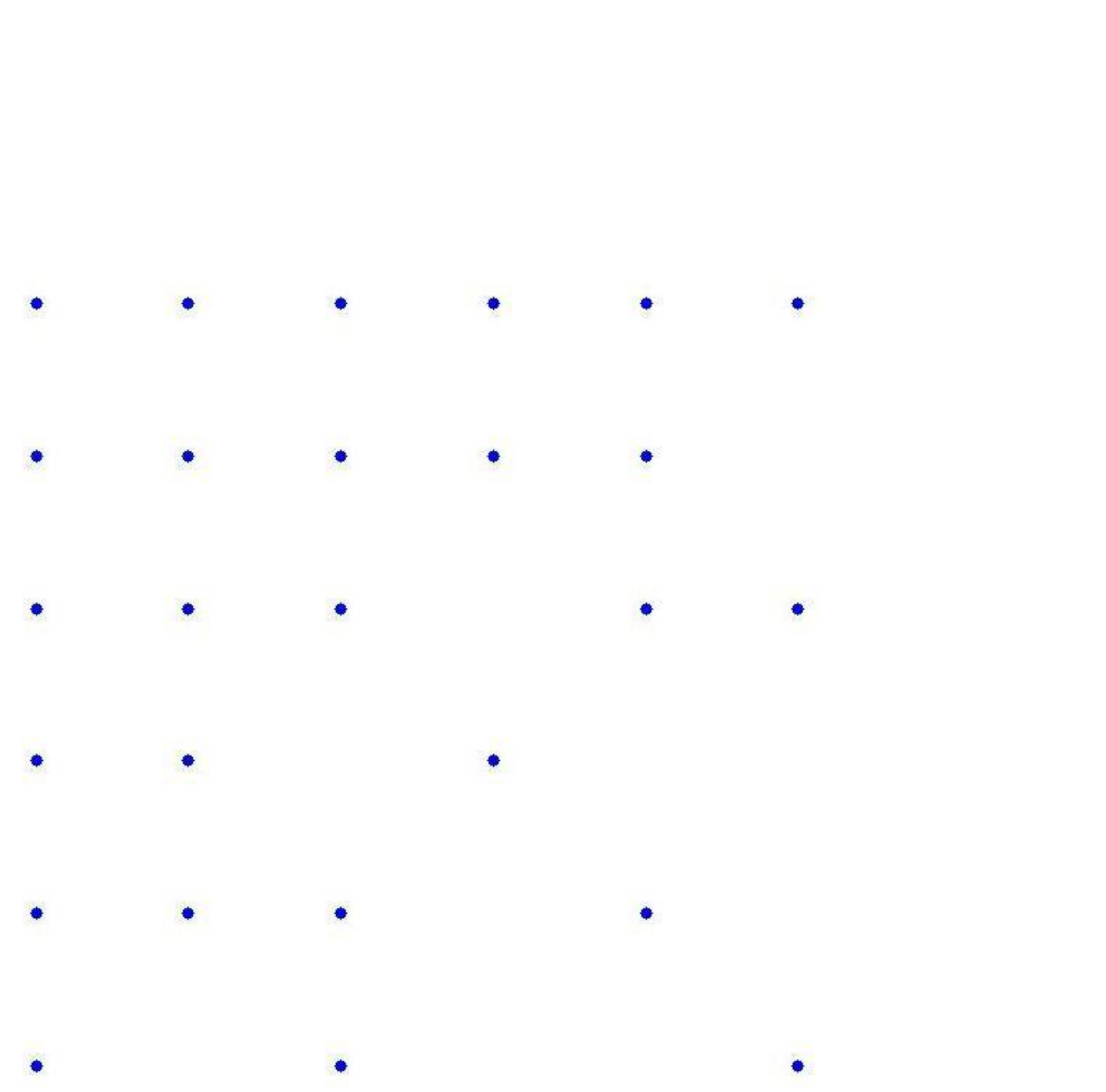}
\hspace{80pt}
\includegraphics[width=142pt,height=116pt]{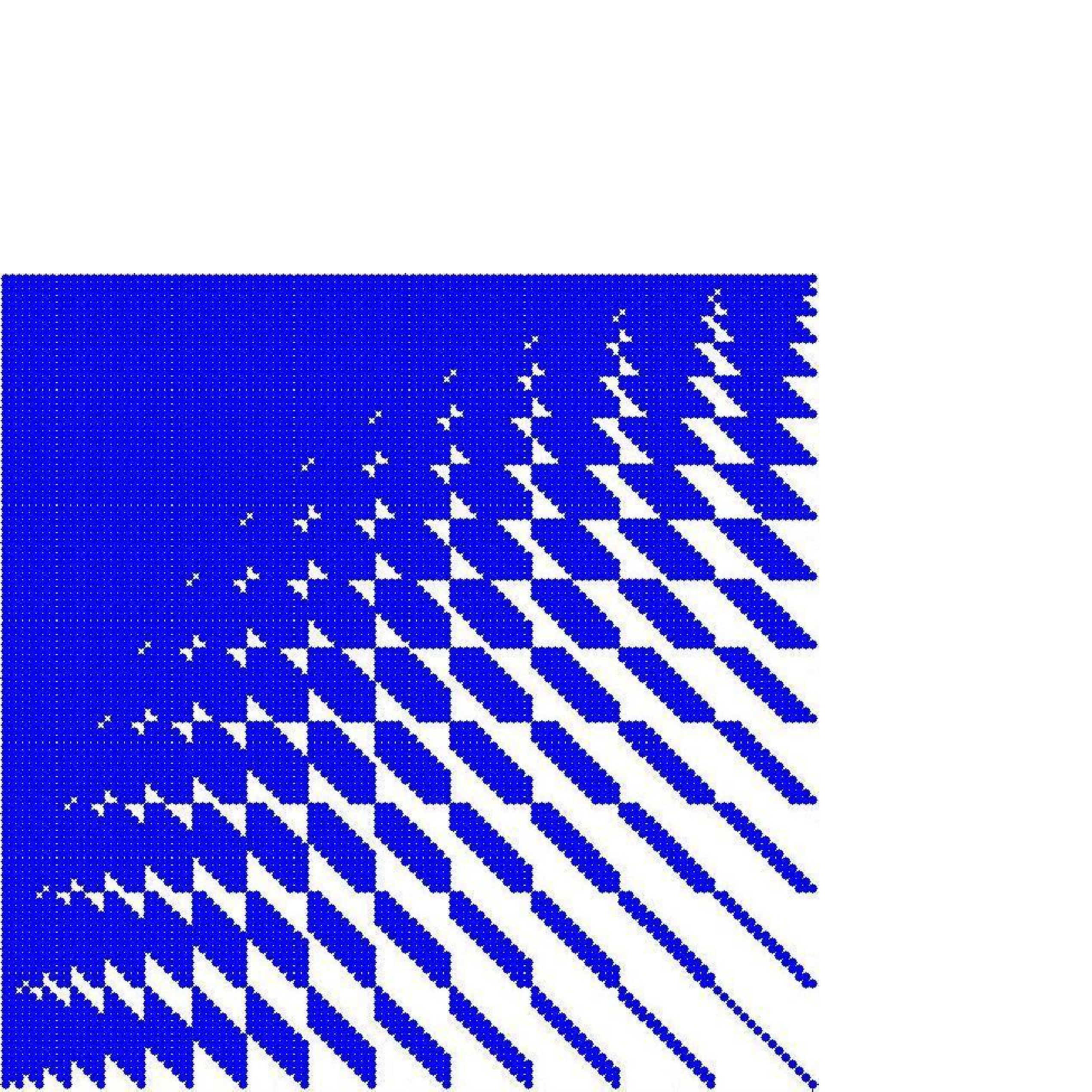}
\caption{Zeros In the Inverse Condition for d=2 and d=14, respectively.  The blue dots indicate non-zero entries.}
\end{figure}

Compare $M_{2}^{-1}$ with Figure \ref{ZIIMM}.

\subsection{First example of parameter collapse to product structure}

Let $\ell$ be a nonnegative integer. Given the bivariate distribution defined on the positive quadrant:

\begin{eqnarray}\label{firstExample}p_{\ell}(x,y)=\frac{1}{C(\ell)} e^{-(x+y)}(x+y)^{\ell}, \end{eqnarray}
where $C(\ell)=\sum_{j=0}^{\ell} \ell! (\ell-j+1)(j+1)$.

Calculate $M_{1}^{-1}$ and set the $(3,2)$ which is the $(2,3)$ entry to zero.  Skipping messy details, this reduces to:\begin{eqnarray*}
0=v(\ell)&\equiv& \left(\int_{0}^{\infty}\int_{0}^{\infty} p_{\ell}(x,y) dxdy\right)\left(\int_{0}^{\infty}\int_{0}^{\infty} xy p_{\ell}(x,y)
dxdy\right)\\
&&-\left(\int_{0}^{\infty}\int_{0}^{\infty} x p_{\ell}(x,y) dxdy\right)^{2}\\
\end{eqnarray*}

It can be shown that $v(\ell)<0$ for $\ell\geq 1$.
Therefore, $\ell=0$ is the only solution which satisfies the zeros in the inverse property for the first truncated moment matrix,
and $p_{0}(x,y)=e^{-x}e^{-y}$ has the product structure indicating independence.  The parameter collapse order of $p_{\ell}(x,y)$ is 1.
Remark that there are a number of similar, parameterized, bivariate, exponential distributions that exhibit parameter collapse.

\subsection{A bivariate exponential non-example}

A multivariate Gamma distribution studied in \cite[Chpt 41\S 3]{4} is:
\[p_{\sigma_{1},\sigma_{2},\rho}(x,y)=\frac{1}{\sigma_{1}\sigma_{2}(1-\rho)} e^{-\frac{1}{1-\rho}\left(\frac{x_{1}}{\sigma_{1}}+\frac{x_{2}}{\sigma_{2}}\right)} I_{0}\left(\frac{2\sqrt{\rho}}{1-\rho}\sqrt{\frac{x_{1}x_{2}}{\sigma_{1}\sigma_{2}}}\right),\]

where $I_{0}$ is the modified Bessel function of first order (with $v=2$).

Solve the collapsing of parameters problems via the Maple code provided in the appendix and the following lines:
\bigskip\hrule\bigskip
\begin{verbatim}
> assume(sigma1>0); assume(sigma2>0); assume(rho <1);
> additionally(rho >0);
> f:= (x,y) -> (1/(sigma1*sigma2*(1-rho)))*
     exp(-(1/(1-rho))*(x/sigma1+y/sigma2))*
     BesselJ(2,2*sqrt(rho*x*y/sigma1/sigma2)/(1-rho));
> equations := zerosInInverse(f,1);
\end{verbatim}
\bigskip\hrule\bigskip
A zero in the $(3,2)$ entry of the first truncated moment matrix implies
\[-(1-\rho)^{2}\left(\rho^{2}+2\rho\left(1- \ln(1+\rho)\right)-2\ln(1+\rho)\right)=0.\]
The left hand side is zero only if $\rho =0$, but then $p_{\sigma_{1},\sigma_{2},0}(x,y)=0$ is not a pdf.
Thus, $p_{\sigma_{1},\sigma_{2},\rho}(x,y)$ never has product structure.

\subsection{An example with compact support}

Consider a probability density function found in \cite[\S 15.4]{5},
\[pdf_{R}(x,y) = \frac{\Gamma((v+2)/2)}{(\pi v)^{1/2} \Gamma(v/2)} \left(1+v^{-1} \left[
                                                                               \begin{array}{cc}
                                                                                 x & y \\
                                                                               \end{array}
                                                                             \right] R \left[
                                                                                         \begin{array}{c}
                                                                                           x \\
                                                                                           y \\
                                                                                         \end{array}
                                                                                       \right]\right),\]
where $R\in SL(2,\mathbb{R}$ and $(x,y)\in B(0,1)\subset \mathbb{R}^{2}$. Ignore the normalizing constant, with $ad-bc=1$,
\[pdf_{R}(x,y)=1+\frac{1}{v}\left(ax^{2}+cxy+bxy+dy^{2}\right).\]

The first order zeros in the inverse condition reduces to $c+b=0$, implying:
\[pdf_{R}(x,y)=1+\frac{1}{v}\left(ax^{2}+dy^{2}\right).\]

The second order zeros in its inverse condition is:
  \[0=3a^2+36a+22da+48+3d^2+36d.\]

From which four cases arise, all of which do not allow for solutions
when the third order \textit{zeros in the inverse condition} is
imposed.  Therefore, parameter collapse does not occur until degree
three.

\subsection{Polynomial and discrete pdfs}

The previous section provided an example of parameter collapse
requiring degree greater than 2.
The computations were easier, as the only product structures allowed
were trivial ones.  Consider measures on $[0,1]\times [0,1]$ of the form
\[\sum_{k=0}^{deg}\sum^{k}_{\ell=0} a_{k,k-\ell} x^{k}y^{k-\ell}.\]

In this case, positivity is required on $[0,1]\times [0,1]$. By integrating, we see:

For simplicity, ignore the above normalization constraint which would yield:
\[1=\int_{0}^{1}\int_{0}^{1}\sum_{k=0}^{deg}\sum^{k}_{\ell=0} a_{k,k-\ell} x^{k}y^{k-\ell}=\sum_{k=0}^{deg}\sum^{k}_{\ell=0} \frac{a_{k,k-\ell}}{(k+1)(k-\ell+1)}\]
\[\Rightarrow a_{00} = 1- \sum_{k=1}^{deg}\sum^{k}_{\ell=0}
\frac{a_{k,k-\ell}}{(k+1)(k-\ell+1)}.\] Consider the measure:
\[pdf^{a_{0,0},a_{1,0},a_{0,1},a_{1,1}}(x,y)=a_{1,1}xy+a_{1,0}x+a_{0,1}y+a_{0,0}\]
The first order \textit{zeros in the inverse condition} requires $a_{11}a_{00}=a_{10}a_{01}$, which then would imply:
\[pdf^{a_{0,0},a_{1,0},a_{0,1},a_{1,1}}(x,y)=\frac{1}{a_{11}}(a_{10}+a_{11}y)(a_{11}x+a_{01})\]
which is of product form.  Parameter collapse occurs in the
first order.  Parameter collapse for the most general polynomial and discrete measures of higher order result in only conditional independence.  One interesting question is, which parameterized higher order polynomials exhibit parameter collapse to full independence? Note that the answer may yield irreducibility criterion.

\section{Conclusion} There are a wealth of examples of parameter collapse due to the zeros in the inverse condition.
One task related to other areas of research is automation of product form recognition for functions.
Another unanswered question related to \cite{6} is, ``Which polynomials may represent pdfs?''
An interesting problem whose solution requires analysis of integral equations is,
``For integer n, what is the general form of pdfs that exhibit parameter collapse order n?''
The latter two questions depend heavily upon the domain of the pdf.

\appendix

\section{Maple Code}
Here is a maple procedure that calculates the truncated moment matrix of order $d$ for the bivariate pdf, $g(x,y)$;
computes the inverse;
and creates a list of equations generated from the zeros in the inverse condition.

\bigskip\hrule\bigskip
\begin{verbatim}
> restart: with(linalg):
> zerosInInverse := proc(g::operator,d::integer)
>
> # set up the structure of the truncated moment matrix
> L:=[]:
> for i from 0 to d do
>    for j from 0 to i do
>       L:=[op(L),[i-j,j]]:
>    end do;
> end do;
>
> #fill truncated moment matrix using appropriate domain for g(x,y)
> B := (i,j) -> int(int(x^i*y^j*g(x,y),x=0..infinity,y=0..infinity);
> A:=matrix(nops(L),nops(L)):
> for i from 1 to nops(L) do
>    for j from 1 to nops(L) do
>       A[i,j] := B(L[i][1]+L[j][1],L[i][2]+L[j][2]);
>    end do;
> end do;
> C:= inverse(A);
>
> #Generate set of equations from zeros in inverse condition
> S:={};
> for i from 1 to nops(L) do
>    for j from i to nops(L) do
>       if (max(L[i][1],L[j][1])+max(L[i][2],L[j][2]) > d)
>       then S:={op(S),C[i,j]=0}; end if;
>    end do;
> end do;
> S
> end proc:
\end{verbatim}
\bigskip\hrule\bigskip

Thus, for a particular pdf, as in example \ref{firstExample}.

\bigskip\hrule\bigskip
\begin{verbatim}
> assume(ell::nonnegint);
> p:= (x,y) -> exp(-(x+y))*(x+y)^ell/(ell+1)!;
> equations = zerosInInverse(p,1);
\end{verbatim}
\bigskip\hrule\bigskip

The next step would be to run the following command.
\bigskip\hrule\bigskip
\begin{verbatim}
> solve(equations);
\end{verbatim}
\bigskip\hrule

\end{document}